\documentclass[twocolumn,11pt]{article}
\usepackage{times}
\usepackage{amsmath,amssymb}
\usepackage{enumerate}
\usepackage[english]{babel}
\usepackage[T1]{fontenc}
\usepackage{color}
\usepackage{eepic,epsfig}
\newtheorem{definition}{Definition}
\newtheorem{remark}[definition]{Remark}
\newtheorem{lemma}[definition]{Lemma}
\newtheorem{theorem}[definition]{Theorem}
\newtheorem{proposition}[definition]{Proposition}
\newtheorem{corollary}[definition]{Corollary}
\newtheorem{exe}[definition]{Example}
\newcommand{\eop}{\hfill $\sqcap\!\!\!\!\sqcup$} 

\newcommand{\N}{\mathbb N}

\def\Nset{\mathbb N}
\def\Rset{\mathbb R}
\def\Cset{\mathbb C}

\def\PX{\mathbb P}

\def\PXv{\PX_{g}}
\def\PXr{\PX_{r}}
\begin{document}
\global\def\refname{{\normalsize \it References:}}
\baselineskip 12.5pt
%
%
%
\title{\LARGE \bf On rank one  perturbations  of Hamiltonian system with periodic coefficients }

\date{}

\author{\hspace*{-10pt}
\begin{minipage}[t]{2.7in} \normalsize \baselineskip 12.5pt
\centerline{MOUHAMADOU  DOSSO}
\centerline{ Universit\'e FHB de Cocody-Abidjan}
\centerline{
  UFR Maths-Info.,}
\centerline{ 22 BP 582 Abidjan  22,}
\centerline{C\^OTE D'IVOIRE}
\centerline{mouhamadou.dosso@univ-fhb.edu.ci}
\end{minipage} \kern 0in
\begin{minipage}[t]{2.7in} \normalsize \baselineskip 12.5pt
\centerline{AROUNA G. Y. TRAORE}
\centerline{Universit\'e FHB de Cocody-Abidjan}
\centerline{UFR  Maths-Info.}
\centerline{22 BP 582 Abidjan  22,}
\centerline{C\^OTE D'IVOIRE}
\centerline{traorearounagogbeyves@yahoo.ca}
\end{minipage}
\\
\\
\begin{minipage}[t]{2.7in} \normalsize \baselineskip 12.5pt
\centerline{JEAN-CLAUDE  KOUA BROU}
\centerline{Universit\'e FHB de Cocody-Abidjan}
\centerline{UFR Maths-Info.}
\centerline{22 BP 582 Abidjan  22,}
\centerline{C\^OTE D'IVOIRE}
\centerline{$k_{-}$brou@hotmail.com}
\end{minipage}
%
%
\\ \\ \hspace*{-10pt}
\begin{minipage}[b]{6.9in} \normalsize
\baselineskip 12.5pt {\it Abstract:}
From a  theory developed  by  C. Mehl, et al., a theory of the rank  one perturbation
of Hamiltonian systems  with periodic coefficients is proposed. It is showed  that
 the rank one perturbation of the fundamental solution
 of Hamiltonian system with periodic coefficients is solution of its  rank one perturbation.
 Some results on the consequences of the strong stability of these types
of systems on their rank one  perturbation is proposed.
Two numerical examples are given to illustrate this theory.
\\[4mm]{\it 2010  Mathematics Subject  Classification :}
15A63, 15A21,  47A55, 93B10, 93C73.
\\ [4mm] {\it Key--Words:}
Eigenvalue, symplectic matrix,  Hamiltonian system, Fundamental solutions,   Perturbation.
\end{minipage}
\vspace{-10pt}}

\maketitle

\thispagestyle{empty} \pagestyle{empty}
%
%
\section{Introduction}
\label{S1} \vspace{-4pt}
Let $J,W\in \Rset^{2N\times 2N}$  be two matrices such  that $J$ is nonsingular and skew-symmetric matrix.
We  say  that the matrix $W$  is $J$-symplectic (or $J$-orthogonal ) if $W^TJW=J$.
 These types of matrices (so-called structured)  usually appear in control theory
  \cite{YS,Lanc_95,hass_99,Brez_02}:
 more precisely in optimal control \cite{hass_99} and in the parametric resonance theory  \cite{god-sad_06,YS}.
   In these areas, these types of matrices are obtained as  solutions of  Hamiltonian
    systems  with periodic coefficients.
    About these systems, that are differential equations  with $P$-periodic coefficients of the below form
\begin{equation}\label{EQS0}
J\dfrac{d X(t)}{dt}= H(t)X(t),\;\; t\in \Rset
\end{equation}
where $J^T=-J$, $(H(t))^T=H(t)=H(t+P)$. The  fundamental solution $X(t)$  of (\ref{EQS0})  i.e. the matrix satisfying
\begin{equation}\label{EQS1}
\left\{\begin{array}{rcl}
J\dfrac{d X(t)}{dt}&=& H(t)X(t),\;\; t\in \Rset_+^*\\
X(0)&=&I_{2N}
\end{array}\right.
\end{equation}
is $J$-symplectic \cite{Dos_06,dos-sad_09,dos-cou_15,YS}  and  satisfies the relationship
 $
 X(t+nP)=X(t)X^n(P)$, $\forall t\in \Rset$   and   $\forall  n\in \N$.
  The solution of the system  evaluated at the period is called the monodromy matrix of the system.
  The eigenvalues of this  monodromy matrix  are  called the multipliers of the system (\ref{EQS1}).
  The following definition permits to classify the multiplies of Hamiltonian system
  \begin{definition}
  \label{defKGL}
  Let $\rho$  be a semi-simple  multiplier  of  (\ref{EQS1})  lying  on  the unit circle. Then
  $\rho$  is called   a multiplier  of  the first (second ) kind  if  the quadratic  form
  $(iJx,x)$  is positive  (negative)  on  the eigenspace  associated  with  $\rho$ .
  When $(Jx,x)=0$, then  $\rho$  is  of mixed  kind.
  \end{definition}
  In  this definition, the notation  $(iJx,x)$  stands for the Euclidean  scalar  product  and $i=\sqrt{-1}$.

  This   other   definition
  proposed  by  S. K. Godunov  \cite{dos-coul_14,dos-coul-sam_13,god_92,god-sad_01}
  gives  another classification of the multipliers of (\ref{EQS1})
   \begin{definition}
      Let $\rho$ be a semi-simple multiplier of (\ref{EQS1}) lying on the unit circle. We say  that  $\rho$  is of the red (green) color or
      in short $r$-multiplier ( $g$-multiplier) if $(S_0 x,x)>0$ ( respectively $(S_0x,x)<0 $ ) on the eigenspace
      associated with $\rho$ where $S_0=(1/2)\left((JX(P))^T+(JX(P))\right)$.
       If $(S_0x,x)=0$, we say that $\rho$  is of  mixed color.
       \label{defcolor}
\end{definition}
From  Definition \ref{defcolor}, Dosso  and Sadkane   obtained  a result of strong stability  of symplectic matrix (see \cite{Dos_06,Dos-Sad_13,dos-coul_14})

\begin{theorem}
  A  symplectic matrix  is strong  stability
if and only if
 \begin{enumerate}
 \item all eigenvalues are  on  the unit circle  ;
  \item the eigenvalues  are  either  red color or green color ;
   \item the subspaces associated of these  deux groups of
   the eigenvalues are well separated.
  \end{enumerate}
   \end{theorem}

   Denote by $\PXr$  and $\PXv$ the spectral projectors  associated with the $r-$eigenvalues  and $r-$eigenvalues of the monodromy matrix
    $X(P)$  of (\ref{EQS1})  and let's  put   $S_r:=\PXr^TS_0\PXr=S_r^T\ge 0$ and
   $S_g:=\PXv^TS_0\PXv=S_g^T\leq 0$  where $S_0=(1/2)\left((X(P)J)+(X(P)J)^T\right)$.
We give the following theorem which gathers all assertions on the strong
 stability of Hamiltonian systems with periodic  coefficients \cite{YS,Dos-Sad_13,Dos_06}.

\begin{theorem}\label{thm11}
The Hamiltonian system (\ref{EQS1})  is strongly  stable if one of the following conditions is satisfied :
\begin{enumerate}
\item If there exists  $\epsilon>0$  such  that  any       Hamiltonian  system  with  $P$-periodic  coefficients  of the form
   $J\dfrac{dx(t)}{dt}=\widetilde{H}(t)x(t)$  and  satisfying
$$
\|H-\widetilde{H}\|\equiv \int_0^T \|H(t)-\widetilde{H}(t)\|dt<\epsilon
$$
is  stable.
  \item The monodromy  matrix $W=X(P)$  of the system (\ref{EQS1})   is strongly  stable
  \item (KGL criterion) the multipliers of the system  (\ref{EQS1})  are either  of  the first kind and either of second kind.  The multipliers of the  first kind  and second  kind of the monodromy matrix should be well
      separated i.e. the quantity
     \begin{align}
    &\delta_{KGL}(X(P))=min\left\{|e^{i\theta_k}-e^{i\theta_l}|\; ; e^{i\theta_k},e^{i\theta_l}\right.\nonumber\\
     &\left.\text{are multipliers  of $(\ref{EQS1})$  of different kinds }\right\}\label{GapKGL}
  \end{align}
  should not  be close to zero.
  \item  the multipliers of the  system  (\ref{EQS1})  are either of the red  color  and either of the green  color.
  The  r-multipliers and g-multipliers of the monodromy matrix should be well separated  i.e. the quantity
  \begin{align}
    &\delta_S(X(P))=min\left\{|e^{i\theta_k}-e^{i\theta_l}|\; ; e^{i\theta_k},e^{i\theta_l}\right.\nonumber\\
     &\left.\text{are $r-$multpliers and $g-$multipliers  of}\;\; (\ref{EQS1})\right\}\label{GapS}
  \end{align}
should not  be close to zero.
\item $S_r\ge 0$, $S_g\leq 0$  and $S_r-S_g>0$
\item $\PXr+\PXv=I$  and $\PXr^TS_0\PXv=0$.
\end{enumerate}
\end{theorem}

The paper is organized  as follows. In  Section \ref{Sec1}  we give some  preliminaries   and   useful  results
to  introduce   the rank  one perturbations of Hamiltonian  systems with periodic coefficients.
More specifically, this section explains what led us to  rank
  one perturbations  of Hamiltonian  system with periodic  coefficients.
Section \ref{Sec2} explains the concept of rank  one perturbation of Hamiltonian systems with coefficients.
In  Section \ref{Sec3} we   analyze the consequences of  strongly stable of Hamiltonian  systems  with
periodic coefficients  on its rank on perturbation.  Section \ref{Sec4}
is devoted to numerical tests.
Finally some concluding remarks are summarized in Section  \ref{Sec5}

Throughout   this paper,  we  denoted the identity and zero matrices  of order $k$   by  $I_k$  and $0_k$  respectively  or
just $I$   and $0$   whenever   it  is clear  from   the context. The 2-norm  of a matrix  $A$   is denoted by $\|A\|$.
The transpose of a matrix  (or vector  ) $U$ is denoted by $U^T$.

\section{ Rank one perturbation of symplectic matrices depending on a parameter}\label{Sec1} \vspace{-4pt}
Let  $W\in \Rset^{2N\times 2N}$  be a $J$-symplectic matrix where $J\in \Rset^{2N\times 2N}$  is skew-symmetric matrix  (i.e. $J^T=-J$)\cite{Meh,Meh_11}.
\begin{definition}
We  call a rank  one perturbation  of the symplectic matrix $W$ any  matrix of the form
$\widetilde{W}=(I+uu^TJ)W$  where $u\in \Rset^{2N}$.
\end{definition}
We  recall in the following proposition some  properties of rank one perturbations of symplectic matrices (see \cite{Yan2}).
\begin{proposition}
Let  $W$  be a $J$-sympectic matrix.
\begin{enumerate}
\item Any  rank one perturbation of $W$ is $J$-symplectic.
\item The invertible  of a rank one perturbation $I+uu^TJ$ of  identity matrix $I$  is  the
matrix $I-uu^TJ$.
\end{enumerate}
\end{proposition}
\noindent
{\bf Proof:} \
See  \cite{Meh,Yan2} for the proof.
\eop

\vspace{2pt}

Let $u$  be a vector of $\Rset^{2N}$. Consider  the following  lemma

\begin{lemma}\label{lem1}
  Consider  the  rank one perturbations    $\widetilde{W}=(I+uu^TJ)W$   of  the  $J$-symplectic matrix   $W$.
  Then
for any   $y\in \Rset^{2N}$, the quadratic  form  $(S_0y,y)$  is  defined  by
\begin{equation}\label{S0t}
(S_0y,y)=(\widetilde{S}_0y,y)-\varphi(y)
\end{equation}
where
\begin{align*}
&\widetilde{S}_0=(1/2)\left((J\widetilde{W})+(J\widetilde{W})^T\right)\\
&\quad\text{and}\\
  &\varphi(y)=(1/2)\left(((Juu^TJW)+(Juu^TJW))^T)y,y\right).
  \end{align*}
\end{lemma}

\noindent
{\bf Proof:} \
Developing  $\widetilde{S}_0$, we have
\begin{align*}
\widetilde{S}_0=& (1/2)\left((JW)+(JW)^T\right)+\\
&(1/2)\left[(Juu^TJW)+(Juu^TJW)^T\right].
\end{align*}
we deduct
\begin{align*}
    &(\widetilde{S}_0y,y)= (S_0y,y)+\\
    &\underbrace{(1/2)\left(\left[(Juu^TJW)+(Juu^TJW)^T\right]y,y\right)}_{\displaystyle \varphi(y)}\\
    &\qquad=(S_0y,y)+\varphi(y)
\end{align*}
\eop

\vspace{2pt}

\begin{corollary}\label{cor1}
Let  $\rho$  be an eigenvalue of  $W$ of modulus 1 and  $y$  an  eigenvector    associated with
$\rho$. Then
   $\rho$  is  an  eigenvalue of red color    (respectively  eigenvalue of green ) if   and only  if
  $\left(\widetilde{S}_0y,y\right)>\varphi(y)$ (respectively   $\left(\widetilde{S}_0y,y\right)<\varphi(y)$).  \\
  However  if   $\left(\widetilde{S}_0y,y\right)=\varphi(y)$, then  $\rho$  is  of mixed color.

\end{corollary}

\noindent
{\bf Proof:} \
According to  lemma \ref{lem1}, we get
$$
(S_0y,y)=(\widetilde{S}_0y,y)-\varphi(y)
$$
From Definition \ref{defcolor}, we have
\begin{itemize}
  \item if $\rho$  is an eigenvalue of red color,
  $$
  (S_0y,y)>0 \Longrightarrow (\widetilde{S}_0y,y)>\varphi(y)\; ;
  $$
  \item if  $\rho$  is an eigenvalue of green color,
  $$
  (S_0y,y)<0 \Longrightarrow (\widetilde{S}_0y,y)<\varphi(y)\; ;
  $$
  \item  if $\rho$  is an eigenvalue of mixed color,
  $$
  (S_0y,y)=0 \Longrightarrow (\widetilde{S}_0y,y)=\varphi(y).
  $$
\end{itemize}
\eop

\vspace{2pt}

 We consider the following   rank one perturbation   of the fundamental  solution $X(t)$
of (\ref{EQS1})

\begin{equation}
\widetilde{X}(t)=(I+uu^T)X(t)
\end{equation}
then  we have the following lemma
\begin{lemma}\label{LemSP}
If $\widetilde{X}(t)$  is a $J$-symplectic matrix  function  such  that   $\text{rank}(\widetilde{X}(t)-X(t))=1,\;\, \forall t>0$,
 then there is a vector function  $u(t)\in \Cset^{2N}\;\;\; \forall t>0$  such  that
 $$
 \widetilde{X}(t)=(I+u(t)u(t)^TJ)X(t),\;\;\forall t\in \Rset
 $$
 Conversely,  for any  vector  $u(t)\in \Cset^{2N}$, the  matrix function $\widetilde{X}(t)$  is $J$-symplectic.
\end{lemma}
\noindent
{\bf Proof:} \
According to  Lemma  7.1 of  \cite[Section 7,p. 18]{Meh}, for all $t>0$, there exists a vector $u(t)\in \Cset^{2N}$  such that
$$\widetilde{X}(t)=(I+u(t)u(t)^TJ)X(t).$$
Moreover, if  $X(t)$  is $J$-symplectic,  $\widetilde{X}(t)$   is also  $J$-symplectic.
\eop

\vspace{2pt}

  This Lemma  leads us to introduce the concept of rank one  perturbation of
   Hamiltonian  systems with periodic   coefficients.

    Now consider, in the follow, that the vector function
   is a vector constant.  We give the following  theorem which extend
    Theorem 7.2  of  \cite[Section 7, p. 19]{Meh}  to matrizant of system (\ref{EQS1}).

\begin{theorem}\label{THR}
Let $J\in \Cset^{2N\times 2N}$ be  skew-symmetric and nonsingular   matrix,
 $(X(t))_{t>0}$ fondamental solution of system (\ref{EQS1})   and $\lambda(t)\in \Cset$
  an  eigenvalue of  $X(t)$ for all $t>0$.
  Assume  that   $X(t)$  has  the   Jordan  canonical  form
  {\tiny $$
  \left(\bigoplus_{j=1}^{l_1} \mathcal{J}_{n_1}(\lambda(t))\right)\oplus
  \left(\bigoplus_{j=1}^{l_2}\mathcal{J}_{n_2}(\lambda(t))\right)
  \oplus\cdots\oplus\left(
  \bigoplus_{j=1}^{l_{m(t)}}\mathcal{J}_{n_{m(t)}}(\lambda(t))\right)\oplus \mathcal{J}(t),
  $$}
where  $n_1>\cdots >n_{m(t)}$  with $m\; :\; \Rset \longrightarrow \Nset^*$
 a  function  of index such that the algebraic multiplicities is $a(t)=l_1n_1+\cdots+l_{m(t)}n_{m(t)}$
   and    $\mathcal{J}(t)$  with  $\sigma(\mathcal{J}(t))\subseteq \Cset\setminus\{\lambda(t)\}$
  contains all  Jordan blocks associated with  eigenvalues  different  from  $\lambda(t)$.
   Furthermore, let  $u\in \Cset^{2N}$   and  $B(t)=uu^TJX(t)$.

\begin{enumerate}
  \item[(1)]  If   $\forall t>0$, $\lambda(t)\not\in \{-1,1\}$, then  generically with respect   to the components  of $u$,  the matrix  $X(t)+B(t)$  has the Jordan canonical  form
     {\tiny  \begin{align}
       \left(\bigoplus_{j=1}^{l_1-1} \mathcal{J}_{n_1}(\lambda(t))\right)&\oplus
  \left(\bigoplus_{j=1}^{l_2}\mathcal{J}_{n_2}(\lambda(t))\right)\oplus\cdots
      \oplus\nonumber\\
      &\left(
  \bigoplus_{j=1}^{l_{m(t)}}\mathcal{J}_{n_{m(t)}}(\lambda(t))\right)\oplus \widetilde{\mathcal{J}}(t),\label{EQS2}
      \end{align}}
      where   $\mathcal{J}(t)$   contains  all the Jordan  blocks of $X(t)+B(t)$   associated
      with  eigenvalues  different  from    $\lambda(t)$.

  \item[(2)]  If  $\exists  t_0 >0$, verifying    $\lambda(t_0)\in \{+1,1\}$, we  have
   \begin{enumerate}
   \item[(2a)]  if  $n_1$  is even,  then generically  with   respect to   the components of $u$,  the matrix\\
     $X(t_0)+B(t_0)$   has  the Jordan  canonical form
     {\tiny  \begin{align*}
      \left(\bigoplus_{j=1}^{l_1-1} \mathcal{J}_{n_1}(\lambda(t_0))\right)&\oplus
  \left(\bigoplus_{j=1}^{l_2}\mathcal{J}_{n_2}(\lambda(t_0))\right)
      \oplus\cdots\oplus\\
      &\left(
  \bigoplus_{j=1}^{l_{m(t)}}\mathcal{J}_{n_{m(t)}}(\lambda(t_0))\right)\oplus \widetilde{\mathcal{J}}(t_0),
      \end{align*}}

      where   $\mathcal{J}(t)$   contains all   the Jordan  of   $X(t)+B(t)$
      associated with eigenvalues  different  from   $\lambda(t)$.
 \item[(2b)]  if  $n_1$  is odd, then  $l_1$  is even   and  generically  with   respect   to   the components of $u$,  the matrix
     $X(t_0)+B(t_0)$   has the Jordan  canonical form

{\tiny \begin{align*}
     \mathcal{J}_{n_1+1}(\lambda(t_0)) &\oplus\left(\bigoplus_{j=1}^{l_1-2} \mathcal{J}_{n_1}(\lambda(t_0))\right)
     \oplus\cdots\oplus\\
     &\left(
  \bigoplus_{j=1}^{l_{m(t)}}\mathcal{J}_{n_{m(t)}}(\lambda(t_0))\right)\oplus \widetilde{\mathcal{J}}(t_0),
      \end{align*}}

      where  $\mathcal{J}(t_0)$   contains all the  blocks of  $X(t_0)+B(t_0)$
       associated with eigenvalues different   from      $\lambda(t_0)$.
       \end{enumerate}
\end{enumerate}
\end{theorem}

\noindent
{\bf Proof:} \
For  all  $t>0$, if $\lambda(t)\not\in \{-1,1\}$,  we have the decomposition  $(\ref{EQS2})$ according to   \cite[Theorem 7.2]{Meh}).
 Other hand, the number of  Jordan  blocks  depend on the variation of $t$.  Thus, this number is a function  of index  $m\; : \; \Rset^+ \longrightarrow \Nset^*$.\\
 For the other two points $(2a)$ and $(2b)$, they show in the same way that items (2) and (3) of  Theorem 7.2  of  \cite[Theorem 7.2]{Meh}) since $X(t_0)+B(t_0)$  is a constant matrix.
\eop

\vspace{2pt}

In reality, the integers  $l_1, ..., l_m(t)$  and indexes $n_1, ..., n_m(t)$  are not constant when t varies.
The number of Jordan blocks and their sizes can varied  in function of   the variation of t.  In Theorem \ref{THR}, we considered
the integers $l_k$  and  $n_k$ constant  $\forall k\in \left\{1,...,m(t)\right\}$  for an index $m(t)$  given. When $t=0$, $\lambda(0)=1$  with $m(0)=2N$ and $l_k=1,\;\;\forall k$. All  Jordan blocks are reduced to $1$.

\section{Rank one perturbations  of    Hamiltonian  system with periodic coefficients}\label{Sec2}\vspace{-4pt}
Let $u$  be a constant vector of $\Rset^{2N}$. $(X(t))_{t\in \ge 0}$ the fundamental  solution of system  \ref{EQS1}.
We  have the following  proposition
\begin{proposition}
Consider the perturbed  Hamiltonian  system
\begin{equation}\label{Eq1}
   J\dfrac{d\widetilde{X}(t)}{dt} =\left[H(t)+E(t)\right]\widetilde{X}(t)
\end{equation}
where
$$
E(t)=(Juu^TH(t))^T+Juu^TH(t)+(uu^TJ)^TH(t)(uu^TJ).
$$
Then   $\widetilde{X}(t)=(I+uu^TJ)X(t)$  is a solution of  system (\ref{Eq1}).
\label{Prop_31}
\end{proposition}

\noindent
{\bf Proof:} \
 By derivation of  $\widetilde{X}(t)$, we  obtain :
\begin{align*}
    J\dfrac{d \widetilde{X}(t)}{dt}=&J(I+uu^TJ)J^{-1}J\dfrac{dX(t)}{dt}\\
    =&J(I+uu^TJ)J^{-1}H(t)X(t),\\
    & \text{according  from  system }\;\; (\ref{EQS1})\\
    =&[H(t)+Juu^TH(t)]X(t)\\
    =&[H(t)+Juu^TH(t)](I+uu^TJ)^{-1}\widetilde{X}(t)\\
    =&[H(t)+Juu^TH(t)](I-uu^TJ)\widetilde{X}(t)\\
    \text{because}\;\;& (I+uu^TJ)^{-1}=(I-uu^TJ)\;\; \; (\text{see \cite{Yan2}})\\
    =&\left[H(t)-H(t)uu^TJ+Juu^TH(t)-\right.\\
    &\left.Juu^TH(t)uu^TJ\right]\widetilde{X}(t)\\
    =&{\tiny \text{$\left[H(t)+\underbrace{(Juu^TH(t))^T+Juu^TH(t)+(uu^TJ)^TH(t)(uu^TJ)}_{\displaystyle E(t)}\right]\widetilde{X}(t)$}}
 \end{align*}

 Hence  the following  perturbed Hamiltonian  equation (\ref{Eq1})
where
\begin{align}
E(t)=&(Juu^TH(t))^T+Juu^TH(t)+\nonumber\\
&(uu^TJ)^TH(t)(uu^TJ)\label{Et}
\end{align}
\eop

\vspace{2pt}

We note that   $E(t)$  is   symmetric and   $P$-periodic  i.e. $E(t)^T=E(t)$  and
  $E(t+P)=E(t)$  for all   $t\in \ge 0$.  The following  corollary gives us  a simplified  form  of    system  (\ref{Eq1})
\begin{corollary}
The  system (\ref{Eq1})  can be  put  at the   form
\begin{equation}\label{Eq2}
\left\{\begin{array}{rcl}
    J\dfrac{d\widetilde{X}(t)}{dt}&=&(I-uu^TJ)^{T}H(t)(I-uu^TJ)\widetilde{X}(t),\\
    & & \\
     \widetilde{X}(0)&=&I+uu^TJ
     \end{array}\right.
\end{equation}
\end{corollary}

\noindent
{\bf Proof:} \
Indeed, developing  \\
$(I-uu^TJ)^{T}H(t)(I-uu^TJ)$  and   we get
\begin{align*}
&(I-uu^TJ)^{T}H(t)(I-uu^TJ)=H(t)+\\
&\underbrace{(J^Tuu^TH(t))^T+J^Tuu^TH(t)+(uu^TJ)^TH(t)(uu^TJ)}_{E(t)}
\end{align*}
and $\widetilde{X}(0)=(I+uu^TJ)X(0)=I+uu^TJ$.
\eop

\vspace{2pt}

We give the following corollary
\begin{corollary}
Let   $(X(t))_{t \ge 0}$  be  the fundamental solution of  system   (\ref{EQS1}).\\
 All solution $\widetilde{X}(t)$   of    perturbed system   (\ref{Eq2})   of  system  (\ref{EQS1}),
 is of the   form  $\widetilde{X}(t)=(I+uu^TJ)X(t)$.
\end{corollary}

\noindent
{\bf Proof:} \
From Proposition \ref{Eq1}   if   $X(t) $  is a  solution  de (\ref{EQS1}), the perturbed    matrix  $W(t)=(I+uu^TJ)X(t)$
is a  solution  of  (\ref{Eq2}). \\
Reciprocally,  for any   solution  $\widetilde{X}(t)$  de $(\ref{Eq2})$,  Let's put
$$
X(t)=(I-uu^TJ)\widetilde{X}(t)
$$
where $u$  is the vector defined in system (\ref{Eq2})

$$\Longrightarrow  \widetilde{X}(t)=(I+uu^TJ)X(t)$$
because $(I+uu^TJ)$   is inverse of the matrix \\$(I-uu^TJ)$ (see \cite{Yan2}).
By replacing this expression  $\widetilde{X}(t)$  in  (\ref{Eq2}),  we  obtain
\begin{align*}
J(I+uu^TJ)\dfrac{d}{dt}X(t)=(I-uu^TJ)^T&H(t)X(t)\\
J(I+uu^TJ)\dfrac{d}{dt}X(t)=(I-uu^TJ)^T&H(t)X(t)\\
 (I-uu^TJ)^{-T} J(I+uu^TJ)\dfrac{d}{dt}X(t)=&H(t)X(t)\\
 \underbrace{(I+uu^TJ)^{T} J(I+uu^TJ)}_{=J}\dfrac{d}{dt}X(t)=&H(t)X(t)\\
 J\dfrac{d}{dt}X(t)=&H(t)X(t)
\end{align*}
and $X(0)=(I-uu^TJ)\widetilde{X}(0)=(I-uu^TJ)(I+uu^TJ)=I$. Consequently, $X(t)$  is   solution of  (\ref{EQS1}).
\eop

\vspace{2pt}

From the foregoing,  we give the following definition :
\begin{definition}
We call rank one perturbations of Hamiltonian system with periodic coefficients, any perturbation of the form
(\ref{Eq2})  of (\ref{EQS1}).
\end{definition}
Consider the following canonical perturbed system  taking  $I_{2N}$  at $t=0$.
\begin{equation}\label{Eq4}
\left\{\begin{array}{rcl}
    J\dfrac{d\widetilde{W}(t)}{dt}&=&(I-uu^TJ)^{T}H(t)(I-uu^TJ)\widetilde{W}(t),\\
    & & \\
     \widetilde{W}(0)&=&I
     \end{array}\right.
\end{equation}

\section{ Consequence of the strong stability  on rank one  perturbations}\label{Sec3}\vspace{-4pt}
We  give the  following  proposition which is  a consequent of Corollary  \ref{cor1}
\begin{proposition}
If a symplectic matrix  $W$  is strongly stable, then there exists a  positif constant $\delta$  such  that
any  vector $u\in \Rset^{2N}$  verifying  $\|uu^TJW\|<\delta$, we have $\left(\widetilde{S}_0y,y\right)\not=\varphi(y)$
for any  eigenvector $y$  of $W$ where $\widetilde{S}_0=(1/2)\left((J\widetilde{W})+(J\widetilde{W})\right)$
with $\widetilde{W}=(I+uu^T)W$.
\end{proposition}

\noindent
{\bf Proof:} \
The strong stability of symplectic matrix $W$ implies that the eigenvalues of $W$  are either of red color  either
of green color i.e. for any eigenvector  $y$  of $W$, we have
 $$(S_0y,y) \neq 0 \Longrightarrow  (\widetilde{S}_0y,y)\neq \varphi(y)$$
using   Corollary  \ref{cor1}.
\eop

\vspace{2pt}

This  following   Proposition   gives us another  consequence of the strong  stability of $W$ under small  perturbation
 that preserve symplecticity.
\begin{proposition}
If a symplectic matrix $W$ is strongly stable, then there exists   a positif  constant $\delta$  such  that
any  vector  $u\in \Rset^{2N}$  verifies  $\|uu^TJW\|<\delta$,  we  have $\widetilde{W}=(I+uu^tJ)W$  is stable.
\end{proposition}

\noindent
{\bf Proof:} \
If $W$  is strongly  stable, then there exists a positif  constant $\delta$  such  that any small perturbation $\widetilde{W}$ of $W$
preserving its symplecticity verifying $\|\widetilde{W}-W\|\leq \delta$, is stable. In particulary, if the perturbation is a  rank one perturbation  with $\widetilde{W}$ of the form $W+uu^TJW$, any vector $u$ verifying $\|uu^TJW\|\leq \delta$ gives $\widetilde{W}$  stable.
\eop

\vspace{2pt}

Hence we have this following result on the strong stability of the
Hamiltonian  systems with periodic coefficients
\begin{proposition}\label{Prop43}
If   Hamiltonian system with periodic coefficients (\ref{EQS1})  is strongly  stable, then     there exists $\varepsilon >0$ such that
for any vector $u$  verifying
  $$\|E(t)\|\leq \varepsilon $$
  where $E(t)$  is defined in (\ref{Et}),
 rank one perturbation  Hamiltonian system (\ref{Eq2}) associated  is stable.
\end{proposition}
 \noindent
{\bf Proof:} \
This proposition is a consequence of Theorem \ref{thm11} using system (\ref{Eq1}) of   Proposition  \ref{Prop_31}.
 \eop

\vspace{2pt}

On the other hand, if the unperturbed system is unstable,  there exits a neighborhood in  which
 any rank one perturbation of  system (\ref{EQS1})   remains unstable.

 \begin{remark}
 The   stability of any small rank one perturbation of a Hamiltonian system with periodic coefficients
  doesn't imply  its strong stability  because we are in a particular case of the perturbation  of the system.
  However it can permit to study the behavior of multipliers of Hamiltonian systems with periodic coefficients.
\end{remark}

\section{ Numerical examples }\label{Sec4}\vspace{-4pt}

\begin{exe}
Consider  the Mathieu  equation
\begin{equation}
J\dfrac{d^2 y(t)}{dt^2}=(a+b\sin(2t))y(t)
\end{equation}
where $a,b\in \Rset$ (see \cite[vol. 2, p. 412]{YS},\cite{dos-coul_14}). Putting
$$
x(t)=\left(\begin{array}{c}
y\\
\dfrac{dy}{dt}\end{array}\right)
,J=\left(\begin{array}{cc}
0 & -1\\
1 & 0\end{array}\right)$$
and
$$ H(t)=\left(\begin{array}{cc}
b\; \sin\; 2t+a & 0\\0 & 1\end{array}\right),
$$
we obtain  the following canonical  Hamiltonian  Equation
\begin{equation}\label{SPM}
J\dfrac{dX(t)}{dt}=H(t)X(t),\quad \forall t\in \Rset, X(0)=I_2,
\end{equation}
where the matrix $H(t)$  is Hamiltonian   and $\pi$-periodic.
Let $u\in \Rset^{2N\times 2N}$  be a random vector in a neighborhood of the zero vector.
Consider  perturbed system  (\ref{Eq2})  of  (\ref{SPM}).
We show that the rank one  perturbation  of the fundamental solution is a solution of perturbed system (\ref{Eq2}).
Consider
$$
\psi(t)=\|\widetilde{X}_1(t)-\widetilde{X}_2(t)\|,\quad \forall t \ge 0
$$
where $\widetilde{X}_1(t)=(I-uu^TJ)X(t)$  and $(\widetilde{X}_2(t))_{t\in \ge 0}$  is  the
 solution of system (\ref{Eq2}).   We show by  numerical examples
  that $(\psi(t))\leq 1.5\; 10^{-14},\;\;\forall  t\in [0,\pi]$.

\begin{itemize}
\item
For $a=7$  and $b=4$, consider the vector  $u=\left(\begin{array}{c}
 0.8913\\
    0.7621
\end{array}\right)$. In Figure \ref{Fig1}, we  consider a  random  vector u  which permits to disrupt
  system (\ref{SPM}) by   the  vectors  $u,10^{-1}u,10^{-2}u$ and $10^{-3}u$.    In  this first figure,  we  note  that $\psi(t)\leq 1.5\;10^{-14}$.  This shows that  $\widetilde{X}_1(t)= \widetilde{X}_2(t)$  for all $t\in [0,\pi]$  i.e. the rank one  perturbation $(\widetilde{X}_1(t))_{t\in [0,\pi]}$  of the fundamental solution of  system (\ref{SPM})  is equal to the solution $(\widetilde{X}_2(t))_{t\in [0,\pi]}$ of rank one perturbation  system (\ref{Eq2}).

\begin{figure}[h!]
\center
\psfig{figure=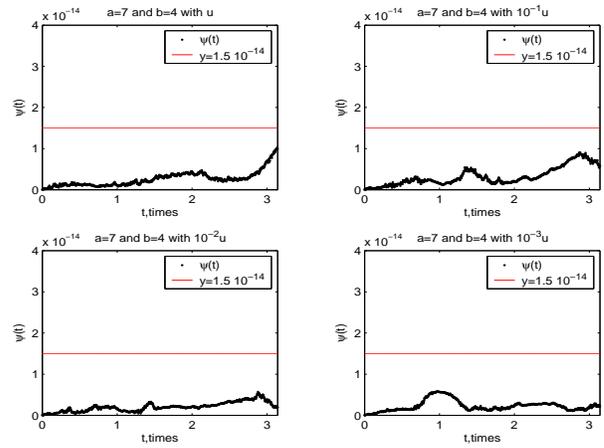,height=60mm,width=80mm}
\caption{Comparison of two  solutions }
\label{Fig1}
\end{figure}

However, unperturbed system  (\ref{SPM}) is strongly  stable.  We  remark that  the rank one perturbed systems (\ref{Eq2})  of (\ref{SPM})   is strongly  stable  when the  vector $u\in \left\{u,10^{-1}u,10^{-2},10^{-3}\right\}$. Therefore   they are stable. This justifies Proposition (\ref{Prop43})

\item  For $a=16.1916618724166685...$  and $b=5$, consider the vector  $u=\left(\begin{array}{c}
0.4565\\
    0.0185
\end{array}\right)$. In this another example illustrated by  Figure \ref{Fig2}, we  consider a  random  vector u  which permits to disrupt
  system (\ref{SPM}) by   the  vectors  $u,10^{-1}u,10^{-2}u$ and $10^{-3}u$. In  figure  \ref{Fig2}, we note  that $\psi(t)\leq 1.5\;10^{-14}$.  This shows that  $\widetilde{X}_1(t)= \widetilde{X}_2(t)$  for all $t\in [0,\pi]$.

\begin{figure}[h!]
\center
\psfig{figure=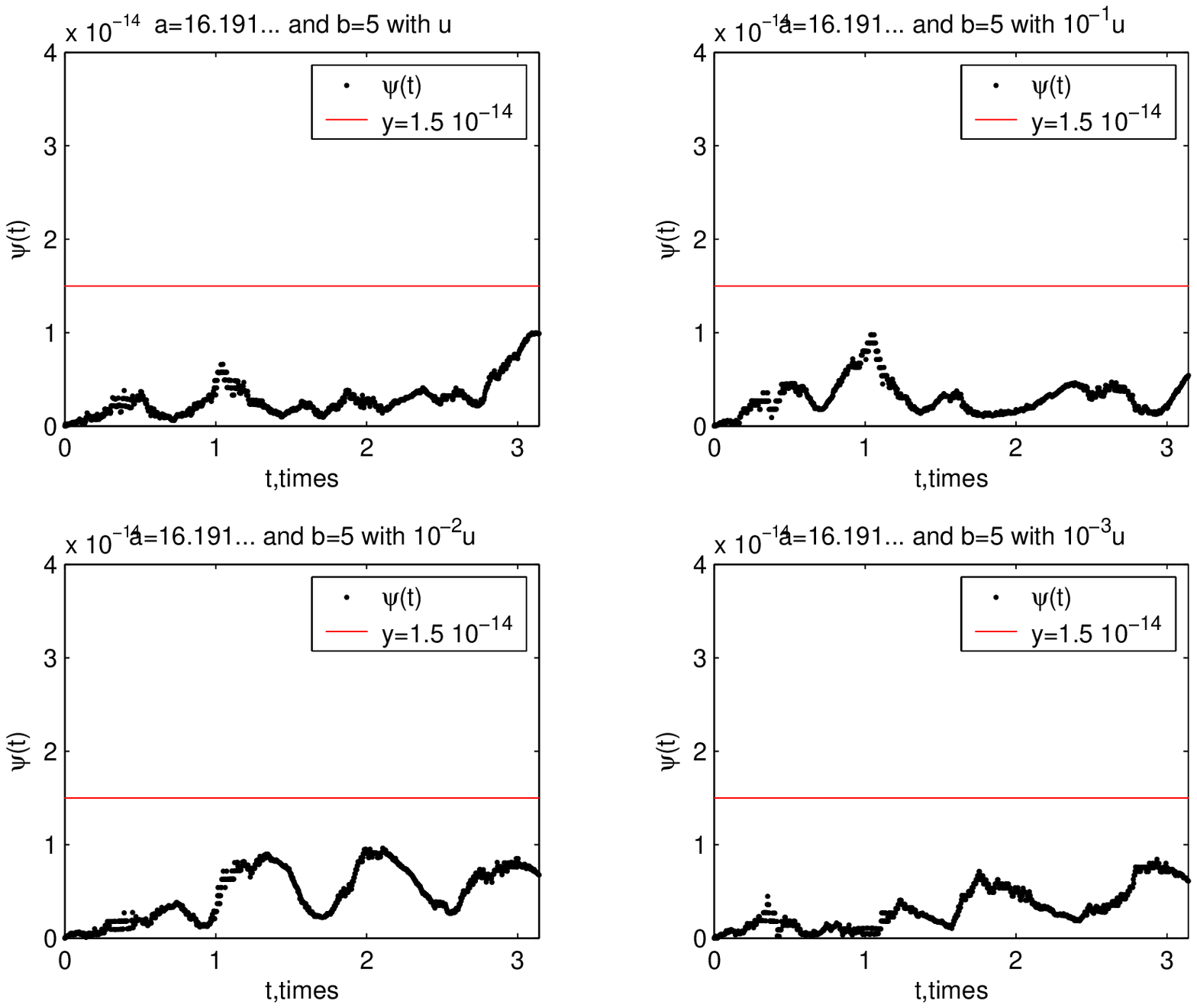,height=60mm,width=80mm}
\caption{Comparison of two  solutions }
\label{Fig2}
\end{figure}

\end{itemize}

In  this example,  the unperturbed system  being unstable,  the rank one perturbation system is unstable when the vector $u\in \left\{u,10^{-1}u,10^{-2}u,10^{-3}\right\}$.  This justifies  the existence of a neighborhood of the unperturbed system
in which any rank one perturbation of the  system is unstable.
\end{exe}

\begin{exe}
Consider the system of differential  equations ( see \cite{god-sad_01}  and \cite[Vil. 2, p. 412]{YS})
\begin{equation}
{\tiny \left\{\begin{array}{l}
q_1\dfrac{d^2\eta_1}{dt^2}+p_1\eta_1+\left[a\eta_1\cos\; 2\gamma t+(b\cos\; 2\gamma t+ c \sin\;2\gamma t)\eta_3\right]=0\\
q_2\dfrac{d^2\eta_2}{dt^2}+p_2\eta_2+g\eta_3\sin\;5\gamma t=0,\\
q_3\dfrac{d^2\eta_3}{dt^2}+p_3\eta_3+\left[(b\cos\;2\gamma t+c\sin\; 2\gamma t)\eta_1+g\eta_2\sin\; 5\gamma t\right]=0,
\end{array}\right.}
\end{equation}
which can be reduced on the following  canonical Hamiltonian system
\begin{equation}\label{SPS}
J\dfrac{dX(t)}{dt}=H(t),\qquad X(0)=I_6
\end{equation}
where
$$
x=\left(\begin{array}{c}
\eta\\
\dfrac{d\eta}{dt}\end{array}\right),\qquad J=\left(\begin{array}{cc}
0_3 & -I_3\\
I_3 & 0_3\end{array}\right),$$
$$ H(t)=\left(\begin{array}{cc}
P(t)  & 0_3\\
0_3 & I_3
\end{array}\right),
$$
with $\eta=\left(\begin{array}{c}
\frac{\eta_1}{\sqrt{q_1}}\\
\dfrac{\eta_2}{\sqrt{q_2}}\\\dfrac{eta_3}{\sqrt{q_3}}\end{array}\right)$
and
{\tiny $$
 P(t)=\left(\begin{array}{ccc}
\dfrac{p_1+a\cos\;2\gamma t}{q_1}  & 0 & \dfrac{b\cos\; 2\gamma\; 2\gamma t+c\sin\;2\gamma t}{\sqrt{q_1q_3}}\\
0  & \dfrac{p_2}{q_2} & \frac{g\sin \: 5\gamma t}{\sqrt{q_2q_3}}\\
\dfrac{b\cos\; 2\gamma\; 2\gamma t+c\sin\;2\gamma t}{\sqrt{q_1q_3}} &
\frac{g\sin \: 5\gamma t}{\sqrt{q_2q_3}} & \dfrac{p_3}{q_3}
\end{array}\right).$$}

Let $u\in \Rset^{2N}$  be a random vector in a neighborhood of the zero vector.
Consider  perturbed system  (\ref{Eq2})  of  (\ref{SPS}).
We show that the rank one  perturbation  of the fundamental solution of  \ref{SPS} is a solution of  its rank  one perturbation  system.
Consider
$$
\psi(t)=\|\widetilde{X}_1(t)-\widetilde{X}_2(t)\|,\forall t\in \Rset
$$
where $\widetilde{X}_1(t)=(I-uu^TJ)X(t)$  and $(\widetilde{X}_2(t))_{t\in \Rset}$  is  the
 solution of the  rank one  perturbation Hamiltonian  system (\ref{Eq2}) of  (\ref{SPS}).  Figures \ref{Fig3}  and \ref{Fig4}   represent the norm of the difference between
 $\widetilde{X}_1$  et $\widetilde{X}_2$.

\begin{itemize}
\item for $\epsilon=15.5$  and $\delta=1$, Let's take\\ $u=\tiny\left(\begin{array}{c}
0.8214\\
    0.4447\\
    0.6154\\
    0.7919\\
    0.9218\\
    0.7382
\end{array}\right)$.  Figure  \ref{Fig3}  is obtained for values of the vector u  taken in $\left\{u,10^{-1}u,10^{-2}u,10^{-3}u\right\}$.  In  figure  \ref{Fig3},  we  note  that $\psi(t)\leq 5\;10^{-13}$.  This shows that  $\widetilde{X}_1(t)= \widetilde{X}_2(t)$  for all $t\in [0,\pi]$  i.e. the rank one  perturbation $(\widetilde{X}_1(t))_{t\in [0,\pi]}$  of the fundamental solution of  system (\ref{SPS})  is equal to the solution $(\widetilde{X}_2(t))_{t\in [0,\pi]}$ of the rank one perturbation  system of (\ref{SPS}).

\begin{figure}[h!]
\center
\psfig{figure=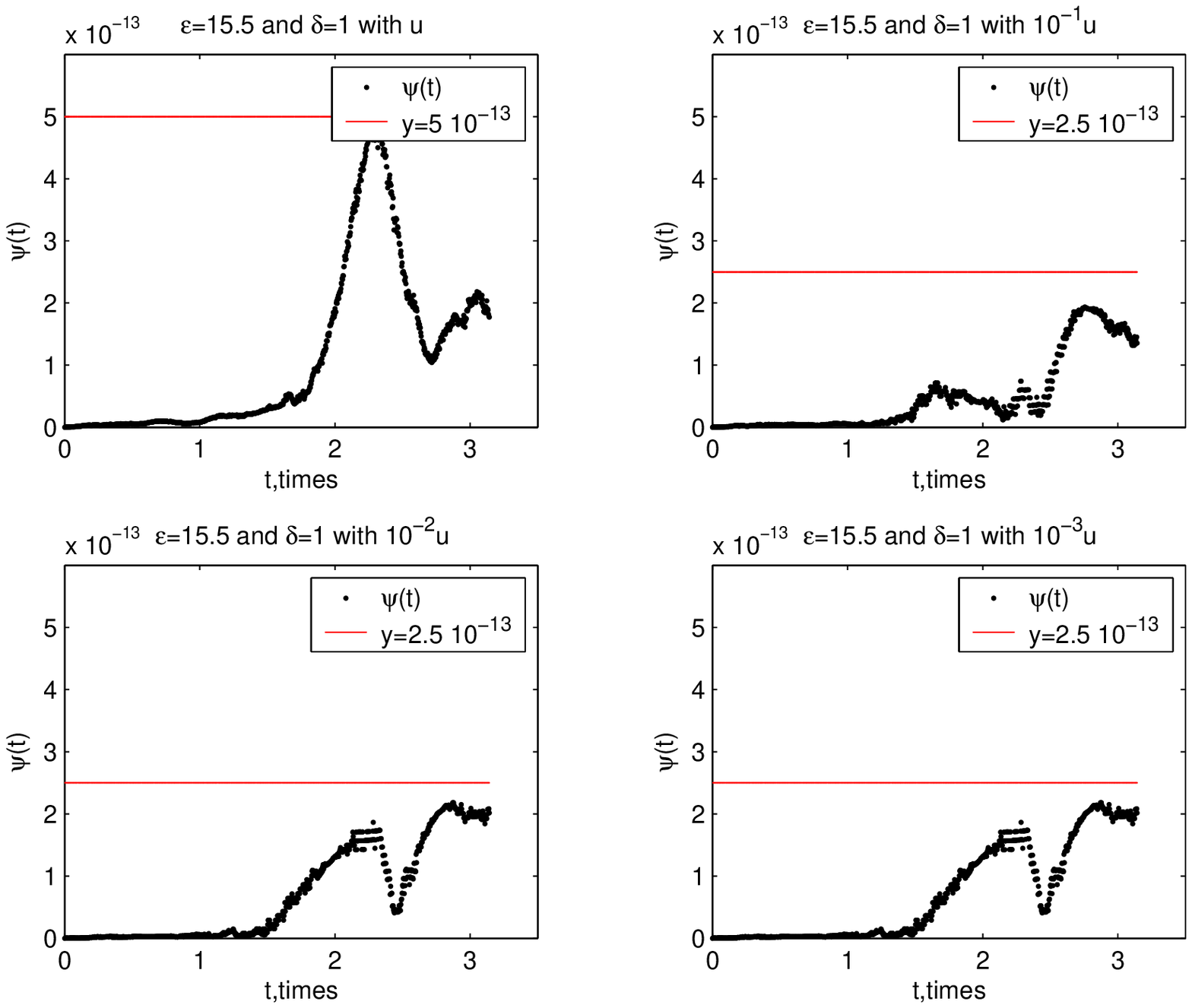,height=60mm,width=80mm}
\caption{Comparison of two solutions }
\label{Fig3}
\end{figure}

However, unperturbed system  (\ref{SPS}) is strongly  stable.  We also note that  the rank one perturbed systems (\ref{Eq2})  of (\ref{SPS})   is strongly  stable  when the  vector $u\in \left\{u,10^{-1}u,10^{-2}u,10^{-3}\right\}$. Therefore   they are stable. This justifies Proposition (\ref{Prop43})

\item  $\epsilon=15$  and $\delta=2$, Let's take $u= \tiny \left(\begin{array}{c}
0.0272\\
    0.3127\\
    0.0129\\
    0.3840\\
    0.6831\\
    0.0928
\end{array}\right)$.  The following figures is obtained for values of the vector u  taken in $\left\{u,10^{-1}u,10^{-2}u,10^{-3}u\right\}$. In  figure  \ref{Fig4}, we  also note   that $\psi(t)\leq 2\;10^{-13}$.  This shows that  $\widetilde{X}_1(t)= \widetilde{X}_2(t)$  for all $t\in [0,\pi]$.
\begin{figure}[h!]
\center
\psfig{figure=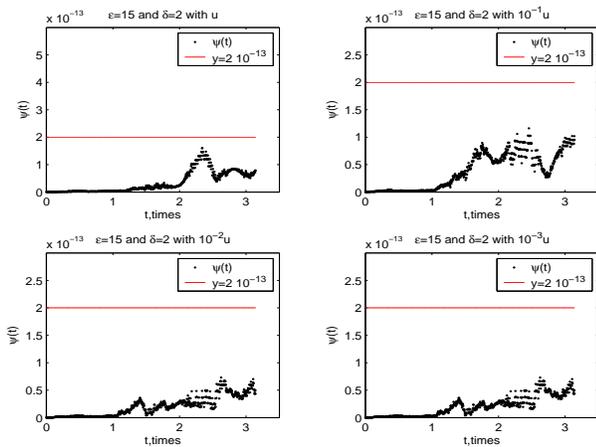,height=60mm,width=80mm}
\caption{ Comparisons of two solutions}
\label{Fig4}
\end{figure}

In  this latter example,  the unperturbed system  is  unstable  and   the rank one perturbation systems remain  unstable when the vector $u\in \left\{u,10^{-1}u,10^{-2}u,10^{-3}\right\}$.  This justifies  the existence of a neighborhood of the unperturbed system
in which any rank one perturbation of the  system is unstable.
\end{itemize}
\end{exe}

\section{Conclusion} \label{Sec5}\vspace{-4pt}
From a theory developed by C. Mehl, et al., on the rank one perturbation of symplectic matrices
 (see \cite{Meh}), we defined the rank one perturbation   of Hamiltonian system of periodic coefficients.
 After an adaptation of some results of \cite{Meh} on  symplectic matrices when they depend on a  time  parameter,
   we show that the rank one perturbation of the  fundamental solution of a Hamiltonian system with periodic coefficients
  is solution of the rank  one perturbation  of the system.  A result of this theory, we give a consequence of the strong stability on a small rank  one perturbation of these Hamiltonian systems.
  Two numerical examples are given to illustrate  this theory.

In future work, we will look how to use the rank one perturbation of
Hamiltonian system with periodic coefficients to analyze  the behavior of
  their  multipliers  and also how this theory can analyze   their  strong stability ?



\begin{thebibliography}{11}
\vspace{-7pt}
\bibitem{Brez_02}
C. Brezinski, Computational  Aspercts  of Linear  Control, {\em Kluwer   Academic
Publishers}, 2002.
\vspace{-7pt}
\bibitem{Dos_06}
M. Dosso, Sur quelques algorithms d'analyse de stabilité forte de matrices
symplectiques, PHD Thesis (September 2006), {\em Université de Bretagne Occidentale.
 Ecole Doctorale SMIS, Laboratoire de Mathématiques, UFR Sciences et Techniques}.

\vspace{-7pt}
\bibitem{dos-cou_15}
M. Dosso, N. Coulibaly, An Analysis  of the Behavior  of Mulpliers
of Hamiltonian  System  with  periodic. {\em Far East  Journal  of Mathematical  Sciences}.
Vol. 99, Number 3, 2016,~301--322.

\vspace{-7pt}
  \bibitem{dos-coul_14}
 M. Dosso, N. Coulibaly, Symplectic matrices  and strong stability
  of Hamiltonian  systems  with  periodic  coefficients. {\em Journal  of Mathematical
   Sciences : Advances  and Applications}.
   Vol.  28, 2014, Pages ~15--38.

\vspace{-7pt}
 \bibitem{dos-coul-sam_13}
 M. Dosso, N. Coulibaly  and L. Samassi, Strong  stability  of symplectic
 matrices  using  a spectral  dichotomy  method. {\em Far  East  Journal  of Applied
 Mathematics}. Vol. 79, Number  2, 2013, pp.~73--110.

\vspace{-7pt}
\bibitem{Dos-Sad_13}
M. Dosso  and M. Sadkane. On the strong stability of symplectic matrices.
  {\em Numerical Linear Algebra with Applications}, 20(2) (2013),~234--249.

\vspace{-7pt}
 \bibitem{dos-sad_09}
M. Dosso, M. Sadkane, A spectral trichotomy method  for symplectic matrices,
{\em Numer Algor.} 52(2009),~187--212


\vspace{-7pt}
 \bibitem{god_92}
S.K. Godunov,   Verification of boundedness for the powers of
symplectic matrices  with the help of averaging, {\em Siber. Math. J.} 33,(1992),
~939--949.

\vspace{-7pt}
\bibitem{god-sad_01}
S.K. Godunov, M. Sadkane,  Numerical determination of a canonical
form of a symplectic matrix, {\em Siberian Math. J.}  42(2001),~629--647.

\vspace{-7pt}
\bibitem{god-sad_06}
S.K. Godunov, M. Sadkane,  Spectral analysis of symplectic matrices
with application to the theory of parametric resonance, {\em SIAM J.
Matrix Anal. Appl.} 28(2006),~1083--1096.

\vspace{-7pt}
\bibitem{hass_99}
B. Hassibi, A.  H.  Sayed, T.  Kailath, Indefinite-Quadratic  Estimation   and Control,
{\em SIAM, Philadelphia}, PA, 1999.

\vspace{-7pt}
\bibitem{Lanc_95}
P. Lancaster, L.  Rodman, Algebraic  Riccati  Equations, {\em Clarendon  Press}, 1995.


\vspace{-7pt}
\bibitem{Meh}
  C. Mehl, V. Mehrmann, A.C.M. Ran, and L. Rodman.  Eigenvalue perturbation
theory under generic rank one perturbations:  Symplectic, orthogonal, and unitary
matrices. {\em BIT}, 54(2014),~219--255.


\vspace{-7pt}
\bibitem{Meh_11}
C. Mehl, V. Mehrmann, A.C.M. Ran and L. Rodman. Eigenvalue perturbation theory of classes
of structured matrices under generic structured rank one perturbations. {\em Linear Algebra Appl.},
435(2011),~687--716.

\vspace{-7pt}
\bibitem{YS}
V.A., Yakubovich,  V.M. Starzhinskii,
\newblock{ Linear differential equations with periodic coefficients},
\newblock Vol. 1 \& 2., Wiley, New York  (1975)
\bibitem{Yan2}
YAN  Qing-you.  The properties  of  a kind  of random  symplectic  matrices.
 Applied  mathematics  and Mechanics. Vol  23, No 5, May  2002.

\end{thebibliography}
\end{document}